\documentstyle[11pt]{article}

\begin{document}

\newcommand{\CC}{\mathop{\rm C\!\!\! I}\nolimits}
\newcommand{\FF}{\mathop{\rm I\! F}\nolimits}
\newcommand{\KK}{\mathop{\rm I\! K}\nolimits}
\newcommand{\LL}{\mathop{\rm I\! L}\nolimits}
\newcommand{\MM}{\mathop{\rm I\! M}\nolimits}
\newcommand{\NN}{\mathop{\rm I\! N}\nolimits}
\newcommand{\PP}{\mathop{\rm I\! P}\nolimits}
\newcommand{\QQ}{\mathop{\rm I\! Q}\nolimits}
\newcommand{\RR}{\mathop{\rm I\! R}\nolimits}
\newcommand{\ZZ}{\mathop{\sl Z\!\!Z}\nolimits}

\newcommand{\pfeil}{\rightarrow}

\newcommand{\kat}{{\cal C}}
\newcommand{\rmat}{{\cal R}}

\newcommand{\horab}{\rule[-1mm]{0pt}{5mm}}
\newcommand{\iso}{\stackrel{\sim}{=}}
\newcommand{\quer}[1]{\overline{#1}}
\newcommand{\schlange}[1]{\widetilde{#1}}
\newcommand{\Nat}{{{\rm Nat}}}
\newcommand{\Mor}{{{\rm Mor}}}
\newcommand{\End}{{{\rm End}}}
\newcommand{\Obj}{{{\rm Obj}}}
\newcommand{\ev}{{{\rm ev}}}
\newcommand{\coev}{{{\rm coev}}}
\newcommand{\id}{{{\rm id}}}
\newcommand{\Id}{{{\rm Id}}}
\newcommand{\Vec}{{{\rm Vec}}}
\newcommand{\Rep}{{{\rm Rep}}}
\newcommand{\Amp}{{{\rm Amp}}}
\newcommand{\RAmp}{{{\rm RAmp}}}

\newcommand{\ZA}{{{\rm ZA}}}
\newcommand{\ZB}{{{\rm ZB}}}
\newcommand{\WA}{{{\rm WA}}}
\newcommand{\WB}{{{\rm WB}}}
\newcommand{\TA}{{{\rm TA}}}
\newcommand{\TB}{{{\rm TB}}}

\newcommand{\akat}{{{\cal A}}}
\newcommand{\bkat}{{{\cal B}}}
\newcommand{\BRib}{{{\rm CylRib}}}
\newcommand{\ARib}{{{\rm Rib}}}
\newcommand{\bnull}{{\overline{b}}}
\newcommand{\dnull}{{\overline{d}}}

\newcommand{\maxleft}{\cap^{-}}
\newcommand{\maxright}{\cap}
\newcommand{\minleft}{\cup^{-}}
\newcommand{\minright}{\cup}
\newcommand{\twist}{\varphi}

\newenvironment{bew}{{\it Proof.}}{\hfill$\Box$}

\newtheorem{bem}{Remark}
\newtheorem{bsp}{Example}	  
\newtheorem{bsps}[bsp]{Examples}
\newtheorem{axiom}{Axiom}
\newtheorem{de}{Definition}
\newtheorem{satz}{Proposition}
\newtheorem{lemma}[satz]{Lemma}
\newtheorem{kor}[satz]{Corollary}
\newtheorem{theo}[satz]{Theorem}

\newcommand{\sbegin}[1]{\small\begin{#1}}
\newcommand{\send}[1]{\end{#1}\normalsize}

\sloppy

\title{Actions of Tensor Categories and Cylinder Braids}
\author{Reinhard H\"aring-Oldenburg\\ 
{\small Mathematisches Institut, Bunsenstr. 3-5,
 37073 G\"ottingen, Germany}\\
{\small email: haering@cfgauss.uni-math.gwdg.de}
}
\date{October 26, 1997}
\maketitle
       
\begin{abstract}
Categorial actions of braided tensor categories are defined and shown
to be the right framework for a discussion of the categorial structure 
related to the group of braids in the cylinder.
A Kauffman polynomial of links in the solid torus is constructed.
\end{abstract}

\section{Introduction}

Braided tensor categories are the great unifying machine of 
braid and link theory. This paper introduces similar notions
for braids in the cylinder and links in the solid torus.

\begin{figure}[ht]
\begin{picture}(150,60)
\put(72,50){\mbox{$\tau_0$}}

\linethickness{0.2mm}
\put(2,50){\mbox{$\cdots$}}
\put(9,56){\mbox{{\small -3}}}
\put(14,56){\mbox{{\small -2}}}
\put(19,56){\mbox{{\small -1}}}
\put(40,56){\mbox{{\small 3}}}
\put(35,56){\mbox{{\small 2}}}
\put(30,56){\mbox{{\small 1}}}

\put(10,45){\line(0,1){10}}
\put(15,45){\line(0,1){10}}
\put(20,55){\line(1,-1){10}}
\put(20,45){\line(1,1){4}}
\put(26,51){\line(1,1){4}}
\put(35,45){\line(0,1){10}}
\put(40,45){\line(0,1){10}}
\put(42,50){\mbox{$\cdots$}}

\linethickness{0.4mm}
\put(90,45){\line(0,1){3}}
\put(90,50){\line(0,1){5}}
\linethickness{0.2mm}
\put(88,51){\oval(4,4)[l]}
\put(88,49){\line(1,0){5}}
\put(93,47){\oval(4,4)[tr]}
\put(93,55){\oval(4,4)[br]}
\put(95,56){\mbox{{\small 1}}}
\put(90,56){\mbox{{\small 0}}}

\put(100,45){\line(0,1){10}}

\put(107,50){\mbox{$\cdots$}}

\put(72,35){\mbox{$\tau_i$}}

\put(0,40){\line(1,-1){10}}
\put(0,30){\line(1,1){4}}
\put(6,36){\line(1,1){4}}
\put(45,40){\line(1,-1){10}}
\put(45,30){\line(1,1){4}}
\put(51,36){\line(1,1){4}}
\put(20,30){\line(0,1){10}}
\put(30,30){\line(0,1){10}}
\put(13,35){\mbox{$\cdots$}}
\put(35,35){\mbox{$\cdots$}}
\linethickness{0.4mm}
\put(90,30){\line(0,1){10}}
\linethickness{0.2mm}
\put(95,30){\line(0,1){10}}
\put(98,35){\mbox{$\cdots$}}
\put(102,40){\line(1,-1){10}}
\put(102,30){\line(1,1){4}}
\put(108,36){\line(1,1){4}}

\put(72,20){\mbox{$e_0$}}

\linethickness{0.2mm}
\put(2,20){\mbox{$\cdots$}}
\put(9,26){\mbox{{\small -3}}}
\put(14,26){\mbox{{\small -2}}}
\put(19,26){\mbox{{\small -1}}}
\put(40,26){\mbox{{\small 3}}}
\put(35,26){\mbox{{\small 2}}}
\put(30,26){\mbox{{\small 1}}}

\put(10,15){\line(0,1){10}}
\put(15,15){\line(0,1){10}}
\put(25,15){\oval(8,8)[t]}
\put(25,25){\oval(8,8)[b]}
\put(35,15){\line(0,1){10}}
\put(40,15){\line(0,1){10}}
\put(42,20){\mbox{$\cdots$}}

\linethickness{0.4mm}
\put(90,15){\line(0,1){10}}
\linethickness{0.2mm}
\put(90,17){\oval(4,4)[tr]}
\put(90,25){\oval(4,4)[br]}
\put(93,26){\mbox{{\small 1}}}
\put(89,26){\mbox{{\small 0}}}

\put(100,15){\line(0,1){10}}

\put(107,20){\mbox{$\cdots$}}

\put(72,5){\mbox{$e_i$}}
\linethickness{0.4mm}
\put(90,0){\line(0,1){10}}
\linethickness{0.2mm}
\put(95,0){\line(0,1){10}}
\put(98,5){\mbox{$\cdots$}}
\put(106,0){\oval(8,8)[t]}
\put(106,10){\oval(8,8)[b]}
\put(20,0){\line(0,1){10}}
\put(30,0){\line(0,1){10}}
\put(12,5){\mbox{$\cdots$}}
\put(32,5){\mbox{$\cdots$}}
\put(5,0){\oval(8,8)[t]}
\put(5,10){\oval(8,8)[b]}
\put(45,0){\oval(8,8)[t]}
\put(45,10){\oval(8,8)[b]}

\end{picture}
\caption{\label{agenerat} 
The graphical interpretation of the generators as symmetric tangles (on the left)
and as cylinder tangles (on the right)	}

\end{figure}
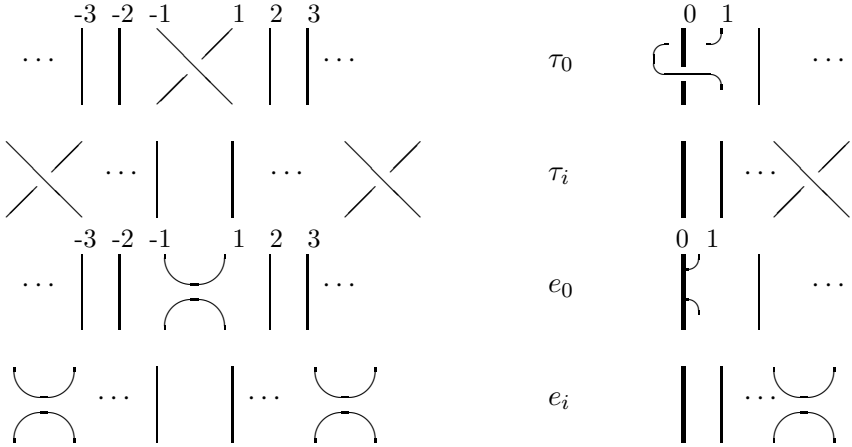

Algebraically, the group of braids in the cylinder appears to be the
braid group related to the Coxeter series B \cite{tD1},\cite{tD2},\cite{lamb}. 
The generators 
$\tau_0,\tau_1,\ldots,\tau_{n-1}$ obey
\begin{eqnarray}
\tau_i\tau_j&=&\tau_j\tau_i\quad\mbox{if}\quad |i-j|>1\\
\tau_i\tau_j\tau_i&=&\tau_i\tau_j\tau_i\quad\mbox{if}\quad i,j\geq1,|i-j|=1
\label{dreizopf}\\
\tau_0\tau_i&=&\tau_i\tau_0\quad\mbox{if}\quad i\geq2\\
\tau_0\tau_1\tau_0\tau_1&=&\tau_1\tau_0\tau_1\tau_0	 \label{vierzopf}
\end{eqnarray}
We denote this group by 
$\ZB_n$. It  may be 
graphically interpreted (cf. figure \ref{agenerat}) as 
symmetric braids or cylinder braids: The symmetric picture
shows it as the group  of 
braids with $2n$ strands (numbered $-n,\ldots,-1,1,\ldots,n$) which are
fixed under a 180 degree rotation about the middle axis. 
In the cylinder picture one adds a single fixed line (indexed $0$)
on the left and
obtains $\ZB_n$ as the group of braids with $n$ strands that may 
surround this fixed line. 
The generators $\tau_i,i\geq0$ are mapped to the 
coresponding diagrams  given in figure \ref{agenerat}.

More generally there are tangles (indicated in figure \ref{agenerat}
by the TLJ tangles $e_i$) of B-type. They are used in the study
of B-type Temperley-Lieb \cite{tD1} and Birman-Wenzl \cite{rhobmw}
algebras.

The need for an extended theory of braided tensor categories
arises because the  braid generator $\tau_0$ cannot be represented by
a morphism in an ordinary braided tensor category. It does not 
satisfy the naturality condition with the  A-type
braiding $\tau_1$.  We account for this fact by separating ordinary morphism
which live in a  braided tensor category from B-type morphisms which 
live in a non-tensor category that is a module over the braided tensor
category. Graphically, the module action is  given by putting 
the ordinary tangle to the right of a cylinder tangle.
This setup has been suggested by tom Dieck \cite{tD1}, \cite{DHO}.

This generality is prompted by the desire to handle morphisms
of the kind of $e_0$ in figure \ref{agenerat}. Restricting 
to tangles that have only braidings around the cylinder one may do with
a somewhat simpler concept introduced in \cite{RHObcat}. 

The primary interest of the present paper 
lies in the formation of concepts. Proofs are rather sketchy,
but may easily be enriched with more details. Physical applictations that
lurk in the background of this work may be found in 
\cite{rhoref}, \cite{rhopott}, \cite{RHObcat}.

  \unitlength1mm


 \unitlength1mm
 \begin{figure}[ht]
\begin{picture}(150,40)

\linethickness{0.4mm}
\put(10,40){\line(0,-1){14}}
\put(10,23){\line(0,-1){17}}
\put(10,3){\line(0,-1){3}}
\linethickness{0.2mm}

\put(15,40){\line(1,-1){10}}
\put(25,40){\line(-1,-1){4}}
\put(15,30){\line(1,1){4}}
\put(15,30){\line(-1,0){4}}

\put(9,27){\oval(6,6)[l]}
\put(9,24){\line(1,0){3}}
\put(12,21){\oval(6,6)[tr]}
\put(25,30){\line(0,-1){10}}

\put(15,21){\line(1,-1){10}}
\put(15,10){\line(1,1){4}}
\put(25,20){\line(-1,-1){4}}

\put(25,11){\line(0,-1){11}}
\put(15,10){\line(-1,0){4}}
\put(9,7){\oval(6,6)[l]}
\put(9,4){\line(1,0){3}}
\put(12,1){\oval(6,6)[tr]}


\put(45,20){\mbox{$=$}}

\linethickness{0.4mm}
\put(60,40){\line(0,-1){5}}
\put(60,32){\line(0,-1){17}}
\put(60,11){\line(0,-1){11}}
\linethickness{0.2mm}

\put(59,36){\oval(6,6)[l]}
\put(59,33){\line(1,0){3}}
\put(62,30){\oval(6,6)[tr]}
\put(75,40){\line(0,-1){10}}

\put(65,30){\line(1,-1){10}}
\put(65,20){\line(1,1){4}}
\put(75,30){\line(-1,-1){4}}

\put(59,16){\oval(6,6)[l]}
\put(59,13){\line(1,0){3}}
\put(62,10){\oval(6,6)[tr]}
\put(75,20){\line(0,-1){10}}

\put(65,10){\line(1,-1){10}}
\put(65,0){\line(1,1){4}}
\put(75,10){\line(-1,-1){4}}

\end{picture}

\caption{\label{vierzopfbild} 
The cylinder interpretation of relation (\ref{vierzopf})}
\end{figure}
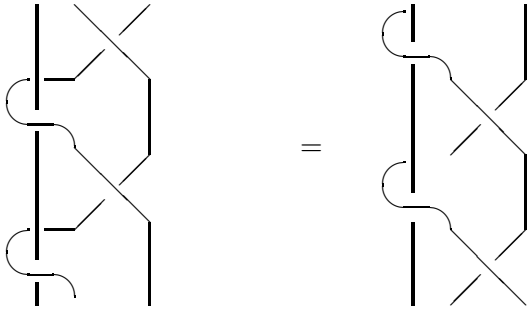

Tammo tom Dieck deserves thanks for discussions which stimulated much 
of the work of this paper.

{\it Preliminaries:}
We use the notation of \cite{ka} for tensor categories.
Expecially we denote by $a_{X,Y,Z}:(X\otimes Y)\otimes Z\rightarrow
X\otimes(Y\otimes Z)$ the associator and by $c_{X,Y}:X\otimes Y\rightarrow
Y\otimes X$ the braiding of a tensor category (resp. braided tensor
category).

\section{Actions of Tensor Categories}

We formalise the notion of a tensor category acting on another category
 in the following way:
\begin{de} Let $\bkat$	be a category and $\akat$ be a tensor category.
We say that $\akat$ acts on $\bkat$ (from the right) if there is a functor
$*:\bkat\times\akat\rightarrow\bkat$ such that the following axioms hold:
\begin{enumerate}
\item The following equation holds whenever both sides are defined:
\begin{equation} (f*g)(f'*g')=(ff')*(gg')
\end{equation}
\item There is a natural isomorphism 
$\lambda\in\Nat(*(\Id\times\otimes),*(*\times\Id))$, 
i.e. $\lambda_{Y,X_1,X_2}:Y*X_1\otimes X_2\rightarrow Y*X_1*X_2$
such that the following pentagon diagram commutes for all
objects $Y\in\Obj(\bkat),X_i\in\Obj(\akat)$:
\begin{equation}
\begin{array}{lcr}
Y*(X_1\otimes X_2)\otimes X_3 & 
\stackrel{\id_Y*a_{X_1,X_2,X_3}}{\longrightarrow} &
Y*X_1\otimes(X_2\otimes X_3) \\
&&\lambda_{Y,X_1,X_2\otimes X_3}\downarrow\\
\downarrow\lambda_{Y,X_1\otimes X_2,X_3}&&Y*X_1*X_2\otimes X_3\\
 &&
\lambda_{Y*X_1,X_2\otimes X_3}\downarrow\\
Y*(X_1\otimes X_2)\otimes X_3 &
\stackrel{\lambda_{Y,X_1,X_2}*\id_{X_3}}{\longrightarrow}&
Y*X_1*X_2*X_3
\end{array}
\end{equation}
\item There is a natural isomorphism $\rho_Y:Y*1\rightarrow Y$
such that 
\begin{equation}\begin{array}{lcr}
Y*1\otimes X & \stackrel{\lambda_{Y,1,X}}{\longrightarrow} & Y*1*X\\
\downarrow \id_Y*l_X&&\rho_Y*\id_X\downarrow\\
Y*X&\stackrel{\id_{Y*X}}{\longrightarrow}&Y*X
\end{array}\end{equation}
Here $1$ denotes the unit object of $\akat$ and 
$l_X:1\otimes X\rightarrow X$ is its compatibility morphism in $\akat$.
\end{enumerate}
The pair $(\bkat,\akat)$ (together with the functor $*$) is called 
an action pair. 
\end{de}

\begin{bsps}\begin{enumerate}
\item  If $F:\akat\rightarrow\bkat$ is a tensor functor between
tensor categories then $\akat$ acts on $\bkat$ by setting
$X*Y:=X\otimes F(Y), \lambda_{X,Y_1,Y_2}:=a^{-1}_{X,F(Y_1),F(Y_2)}$.
As a special case any tensor category acts on itself.
\item Let $\akat$ be the category of bimodules over some ring $R$.
This tensor category acts on the category $\bkat$ of $R$ right 
modules in the obvious way. This example is a special case 
of the former where the functor $F$ is the forgetful functor from the 
category of bimodules to the category of right modules.
\item Let $\akat$ be a group considered as a tensor category, i.e.
the objects are the group elements, tensor product is group 
multiplication. The endomorphism space of an object is some
unital ring $R$ while only one morphism $0\in R$ 
exists between different objects.
Assume that this group acts on a space $\bkat$ which we consider as a
category in a similar way. Then $\akat$ acts on $\bkat$
in the sense of the above definition. This action is strict
according to the definition given below.
\end{enumerate}
\end{bsps}
Further examples will be given later on. 

\begin{de} The action pair $(\bkat,\akat)$ 
is called strict if $\akat$ is a strict tensor category
and  one has
$Y*X_1*X_2=Y*X_1\otimes X_2$, $\lambda_{Y,X_1,X_2}=\id_{Y*X_1*X_2}$
and $\rho_Y=\id_Y$.
\end{de}

\begin{de}
 Let ($\bkat,\akat)$ and $(\bkat',\akat')$ be two
 action pairs. 
A functor between $(\bkat,\akat)$ and $(\bkat',\akat')$ consists
of:
\begin{enumerate}
\item A functor $F_\bkat:\bkat\rightarrow\bkat'$
\item A tensor functor $F_\akat:\akat\rightarrow\akat'$ with
functorial morphisms $\varphi_0,\varphi_2$ defined as in \cite{ka}[XI.4.1].
\item Natural isomorphisms $\omega_{X,Y}:F_\bkat(X*Y)\rightarrow
F_\bkat(X)*F_\akat(Y)$ such that  the following diagram commutes
\begin{equation} \label{funkhexagon}
\begin{array}{lcr}
F_\bkat(Y*X_1\otimes X_2)&\stackrel{\lambda_{Y,X_1,X_2}}{\longrightarrow} &
F_\bkat(Y*X_1*X_2)\\
\downarrow\omega_{Y,X_1\otimes X_2}&&\omega_{Y*X_1,X_2}\downarrow\\
F_\bkat(Y)*F_\akat(X_1\otimes X_2)&&F_\bkat(Y*X_1)*F_\akat(X_2)\\
\downarrow\id*\varphi_2(X_1,X_2)^{-1}&&\omega_{Y,X_1}*\id_{F_\akat(X_2)}\downarrow\\
F_\bkat(Y)*F_{\akat}(X_1)\otimes F_\akat(X_2)&
\stackrel{\lambda'_{F_\bkat(Y),F_\akat(X_1),F_\akat(X_2)}}{\longrightarrow} &
F_\bkat(Y)*F_\akat(X_1)*F_\akat(X_2)
\end{array}
\end{equation}
\item The following diagram commutes
\begin{equation}\begin{array}{lcr}
F_\bkat(Y*1)&\stackrel{F_\bkat(\rho_Y)}{\longrightarrow}&F_\bkat(Y)\\
\downarrow\omega_{Y,1}&&\rho'_{F_\bkat(Y)}\uparrow\\
F_\bkat(Y)*F_\akat(1)&\stackrel{\id*\varphi_0^{-1}}{\longrightarrow} &
F_\bkat(Y)*1
\end{array}\end{equation}
\end{enumerate}
\end{de}

Tensor categories can always be turned into strict ones	 by 
a procedure due to MacLane.
A similar result holds in our situation:
\begin{satz} Every action pair $(\bkat,\akat)$ is equivalent to a strict
action pair  $(\bkat^{\rm str},\akat^{\rm str})$.
\end{satz}
\begin{bew}
The proof is a variation of the proof of MacLanes's theorem.
Hence we restrict ourselves to a sketchy description.

The objects of  $\bkat^{\rm str}$ are sequences of one 
object of $\bkat$ and arbitrary many objects from $\akat$, i.e. 
\begin{eqnarray*}
\Obj(\bkat^{\rm str})&:=&\{(Y,X_1,\ldots,X_k)\mid Y\in\Obj(\bkat),
X_i\in\Obj(\akat),k\in\NN_0\}\\
\Obj(\akat^{\rm str})&:=&\{(X_1,\ldots,X_k)
\mid X_i\in\Obj(\akat),k\in\NN_0\} \end{eqnarray*}
The equivalence functor is defined on objects by 
\begin{eqnarray*}
F_\akat&:&\akat^{\rm str}\rightarrow\akat,(X_1,\ldots,X_k)\mapsto
X_1\otimes(X_2\otimes(\cdots))\\ 
F_\bkat&:&\bkat^{\rm str}\rightarrow\bkat,(Y,X_1,\ldots,X_k)\mapsto
Y*X_1*\cdots*X_k\end{eqnarray*}
Morphism spaces are defined by 
\begin{eqnarray*}
\Mor_\akat^{\rm str}(S_1,S_2)&:=&\Mor_\akat(F_\akat(S_1),F_\akat(S_2))\\
 \Mor_\bkat^{\rm str}(S_1,S_2)&:=&\Mor_\bkat(F_\bkat(S_1),F_\bkat(S_2))\end{eqnarray*}
 The functors $F_\bkat,F_\akat$ are essentially faithful and fully faithful.
Hence, they are  equivalences of categories.
Their right inverses are defined by $Y\mapsto(Y)$. 

Tensor product and action are defined by joining sequences. 
It remains to exhibit the natural isomorphism 
$\omega_{S,S'}:F_\bkat(S*S')\rightarrow F_\bkat(S)*F_\akat(S')$.
Its definition is recursive on the length of $S'$. One sets
\[\omega_{S,()}:=\rho_{F_\bkat(S)}^{-1},\omega_{S,(X)}:=\id,
\omega_{S,(X)\otimes S'}:=\lambda^{-1}_{F_\bkat(S),X,F_\akat(S')}
\omega_{S*(X),S'}\] 
The key lemma to establish (\ref{funkhexagon}) is 
\begin{lemma}
$\lambda_{F_\bkat(S),F_\akat(S'),F_\akat(S'')}
(\id_{F_\bkat(S)}*\varphi_2(S',S'')^{-1})\omega_{S,S'\otimes S''}
=(\omega_{S,S'}*\id_{F_\akat(S'')})\omega_{S*S',S''}
$\end{lemma}
It is shown by induction on the length of $S'$.
\end{bew}

The strictification of action pairs simplifies considerably the
task of specifying them by generators and relations in a
fashion  similar to the presentation of braided tensor 
categories given in \cite{ka}[XII.1].
One starts with a strict action pair $(\bkat,\akat)$ and singles out a 
set ${\cal F}_\bkat$ of morphisms from $\bkat$.
They are used to build formal words defined recursively by their length:
Words of length $1$ are $[f]$ where $f\in{\cal F}$ and 
$[\id_Y], Y\in\Obj(\bkat)$. If $a,b$ are words of length $\leq n$
and $g$ is a morphism from $\akat$
then $a*g$ and $ab$ are words of length $n+1$. To every word
a morphism of $\bkat$ is associated by the rules
$\overline{[f]}:=f,\overline{a*g}:=\overline{a}*g,
\overline{ab}:=\overline{a}\circ\overline{b}$.
The set of sub-words of a word is also defined recursively by 
${\rm sub}([f]):=\{[f]\}, {\rm sub}(a*b):=\{b\}\cup{\rm sub}(a),
{\rm sub}(ab):={\rm sub}(a)\cup{\rm sub}(b)$.
Two words $a,$ are said to be equivalent $a\sim b$ iff
there exists a sequence of words $a_i$ with $a_0=a,a_k=b$ and
$a_{i+1}$ is obtained from $a_i$ by one of the following transformations 
applied to a sub-word: 
$(ab)c\sim a(bc),[\id]a\sim a,a[\id]\sim a,a*\id_1\sim a,
[\id_{Y*X}]\sim[\id_Y]*\id_X, a*gg'\sim a*g*g',
(a*g)(a'*g')\sim(aa')*(gg')$. From this one concludes
that $(a*\id_{b(g)})([\id_{s(\overline{a})}]*g)\sim
([\id_{b(\overline{a})}]*g)(a*\id_{s(g)})$ and
$(a_1*\id)\cdots(a_k*\id)\sim (a_1\cdots a_k)*\id$.
A simple inductive proof shows that any word is equivalent to one
of the form $h_1\cdots h_m$ where each $h_i$ is of the form
$[f]*\id_{X}$ with $f\in{\cal F}$ or of the form $[\id_X]*g$.

The free action pair generated by ${\cal F}$ is the pair
$({\cal M}({\cal F}),\akat)$ where ${\cal M}({\cal F})$ has
the same objects as $\bkat$ but its morphism space is 
the set of equivalence classes of words. 

Further relations ${\cal R}=\{(r_i,r'_i)\mid i=1..k\}$ can be
used to define another equivalence relation $a\sim_{\cal R}b$ on
words  where one may also replace a sub-word $r_i$ by $r'_i$
or vice versa. One then says that the action pair $(\bkat,\akat$)
is generated by ${\cal F}$ with relations ${\cal R}$ if
every morphism of $\bkat$ can be obtained 
as $\overline{a}$ from a word and one has $a\sim_{\cal R}b
\Leftrightarrow \overline{a}=\overline{b}$.

\section{Cylinder twists}

This section introduces the cylinder braid morphism.

\begin{de} A strict action pair $(\bkat,\akat)$ is said to be 
cylinder braided if: 
\begin{enumerate}
\item $\Obj(\bkat)=\Obj(\akat)$ and $1*X=X$
\item  $\akat$ is a braided tensor category with
braid isomorphisms $c_{X,Y}\in\Mor_\akat(X\otimes Y,Y\otimes X)$
\item  For every object there exists an isomorphism 
$t_X\in\Mor_\bkat(X,X)$ such that
\begin{eqnarray}
c_{Y,X}(t_Y\otimes\id_X)c_{X,Y}(t_X\otimes\id_Y)&=&
(t_X\otimes\id_Y)c_{Y,X}(t_Y\otimes\id_X)c_{X,Y}=t_{X\otimes Y}
\label{ctwist}\\
ft_X&=&t_Yf\quad\forall f\in\Mor_\akat(X,Y)
\end{eqnarray}
\item The following equations should hold if 
$\akat$ is equipped with a duality.
\begin{eqnarray}
(t_{X}\otimes\id_{X^\ast})b_X&=&
c^{-1}_{X,X^\ast}(t^{-1}_{X^\ast}\otimes\id_X)
c^{-1}_{X^\ast,X}b_X\label{btrel}\\
d_X(t_X^{\ast-1}\otimes\id_{X})&=&
d_Xc_{X,X^\ast}(t_X\otimes\id_{X^\ast})c_{X^\ast,X}	\label{dtrel}
\end{eqnarray}
\end{enumerate}
$t$ is called the cylinder twist. For the sake of brevity we call
$(\bkat,\akat)$ (or even $\bkat$) a cylinder braided tensor category CBTC.
\end{de}
The requirements of strictness and those of point 1 of the definition 
imply that $X*Y=1*X*Y=1*X\otimes Y=X\otimes Y$.  
Note also that in the light of (\ref{ctwist}) relations
(\ref{btrel}), (\ref{dtrel}) may be rewritten as
$t_{X\otimes X^\ast}b_X=b_X$ and
$d_Xt_{X^\ast\otimes X}=d_X$.

\begin{bem}
\begin{enumerate}
\item The space $\End_\bkat(X^{\otimes n})$ carries a representation of the 
braid group $\ZB_n$.
\item Assume there are $m$ distinct morphisms $t^{(1)},\ldots,t^{(m)}$ such that
each product of  pairwise different $t^{(i)}$  
makes the action pair $(\bkat,\akat)$ cylinder braided. 
Then one has a representation of the braid group of the handlebody
of genus $m$ \cite{sos}.
\item Our  action pairs  are defined by a right action of $\akat$.
Similarly one can consider left actions. Suppose $\akat$ to act on $\bkat$
from the right and on $\bkat´$ from the left. If both of these actions
are cylinder braided then one has a tensor representation of the 
braid group of the affine Coxeter diagram
$\bullet=\bullet-\bullet-\cdots-\bullet-\bullet=\bullet$.
\end{enumerate}
\end{bem}

The fundamental geometric example of tangles in the cylinder
will be described in the next section. Here we restrict ourselves
to some simple examples.
\begin{bsps}\begin{enumerate}
\item A ribbon category $\akat$ acting on itself 
is trivially a cylinder braided pair where the cylinder twist
is given by the ribbon twist $t_X=\theta_X$.
\item An abelian group $G$ together with bilinear pairing 
$c:G\times G\rightarrow K^\ast$ with values
in the group of units of a commutative unital ring $K$ may be
viewed as a braided tensor category $\akat$ as in \cite{tu}[p. 29]. 
The pair $(\akat,\akat)$ is then cylinder braided if there is a map
$t:G\rightarrow K^\ast$ such that $t(gg')=c(g,g')c(g',g)t(g)t(g')$. 
In the symmetric case $t$ is simply a group character.
\item Let $\akat$ be a tensor category and $X\in\Obj(\akat)$ any object.
This category acts on $\bkat$ which has the same objects and
morphisms $\Mor_\bkat(X_1,X_2):=\Mor(X\otimes X_1,X\otimes X_2)$. 
The action 	is given by the monoidal product and cylinder twist is
$t_Y:=c_{Y,X}c_{X,Y}$.
\end{enumerate}\end{bsps}

Further examples are provided by the Coxeter-B braided tensor categories 
studied in \cite{RHObcat}. That paper contains a discussion of
cylinder braid structures on Hopf algebras and Tannaka-Krein duality.

\section{Cylinder Ribbon Tangles}

The fundamental example of a cylinder braided action pair is 
the pair $(\BRib,\ARib)$. $\ARib$ is  Turaev's
category of ribbon tangles and $\BRib$ is  the category of 
cylinder ribbon tangles which is  defined just like $\ARib$ but with
the restriction that the tangles extend only in the space
$(\RR^2-(0,0))\times{[}0,1]$.
The action of a tangle $f$ from $\ARib$ on a tangle $g$ from $\BRib$ is
given by putting $f$ to the right of $g$. The category of $\akat$ coloured
cylinder ribbon tangles $\BRib_\akat$ parallels the category
$\ARib_\akat$.

We use Turaev's notation for the generators of
$\ARib$. The basic generators of $\BRib$ are 
$\tau^{\downarrow\pm}\tau^{\uparrow\pm}$. They are oriented versions 
of $\tau_0$ given in  figure \ref{agenerat} and its inverse. The arrow
indicates the orientation. 	The lines
are meant to represent ribbons that with the framing  oriented towards 
the axis.

\begin{satz}
The following list of relations holds in $\BRib$:
\begin{eqnarray}
\tau^{\downarrow+}&=&{\tau^{\downarrow-}}^{-1}\label{satz51}\\
\tau^{\uparrow+}&=&{\tau^{\uparrow-}}^{-1}\label{satz52}\\
\tau^{\downarrow-}&=&(\maxright\otimes\downarrow)
({\twist'}^\uparrow\otimes\downarrow\otimes\downarrow)
(\tau^{\uparrow+}\otimes X^-)(\minleft\otimes\downarrow)	\label{satz53}\\
\tau^{\uparrow-}&=&(\maxleft\otimes\uparrow)
(\twist'\otimes\uparrow\otimes\uparrow)
(\tau^{\downarrow+}\otimes T^-)(\minright\otimes\uparrow)	\label{satz54}\\
(\tau^{\downarrow+}\otimes\downarrow)X^+ 
(\tau^{\downarrow+}\otimes\downarrow)X^+&=& 
X^+(\tau^{\downarrow+}\otimes\downarrow) X^+
(\tau^{\downarrow+}\otimes\downarrow)\label{satz55}\\
(\tau^{\uparrow+}\otimes\uparrow)T^+ 
(\tau^{\uparrow+}\otimes\uparrow)T^+&=& 
T^+(\tau^{\uparrow+}\otimes\uparrow) T^+
(\tau^{\uparrow+}\otimes\uparrow)\label{satz56}\\
(\tau^{\downarrow+}\otimes\uparrow)Y^- 
(\tau^{\uparrow+}\otimes\downarrow)Z^-&=& 
Y^-(\tau^{\uparrow+}\otimes\downarrow) Z^-
(\tau^{\downarrow+}\otimes\uparrow)\label{satz57}\\
(\tau^{\uparrow+}\otimes\downarrow)Z^- 
(\tau^{\downarrow+}\otimes\uparrow)Y^-&=& 
Z^-(\tau^{\downarrow+}\otimes\uparrow) Y^-
(\tau^{\uparrow+}\otimes\downarrow)\label{satz58}\\
(\tau^{\downarrow+}\otimes\uparrow)\minright&=&
Y^+(\tau^{\uparrow-}\twist^\uparrow\otimes\downarrow)
\minleft\label{satz59}\\
\maxright(\tau^{\uparrow-}\otimes\downarrow)&=&
\maxleft(\tau^{\downarrow+}\twist'\otimes\downarrow)Y^-\label{satz510}\\
(\uparrow\otimes\twist)\minleft&=&(\tau^{\uparrow+}\otimes\downarrow)
Z^-(\tau^{\downarrow+}\otimes\uparrow)\minright\label{satz511}\\
\maxright(\uparrow\otimes\twist)&=&
\maxleft(\tau^{\downarrow+}\otimes\uparrow)Y^-
(\tau^{\uparrow+}\otimes\downarrow)\label{satz512} \\
(\downarrow\otimes\twist^\uparrow)\minright&=&
(\tau^{\downarrow+}\otimes\uparrow)
Y^-(\tau^{\uparrow+}\otimes\downarrow)\minleft\label{satz513}\\
\maxleft(\downarrow\otimes\twist^\uparrow)&=&
\maxright(\tau^{\uparrow+}\otimes\downarrow)Z^-
(\tau^{\downarrow+}\otimes\uparrow)\label{satz514} 
\end{eqnarray}
\end{satz}

\unitlength1mm
 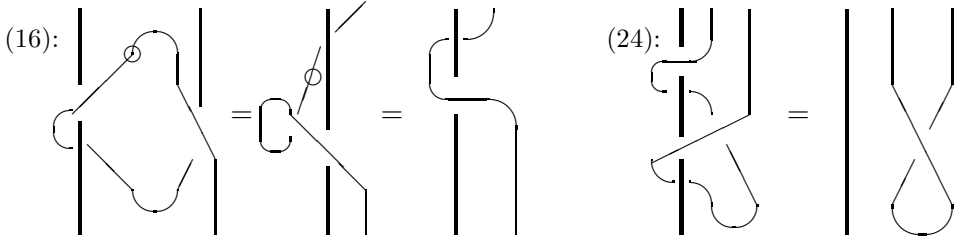
\begin{figure}[ht]
\begin{picture}(150,30)

\put(0,25){\mbox{{\small (\ref{satz53}):}}}
\linethickness{0.4mm}
\put(10,0){\line(0,1){15}}
\put(10,20){\line(0,1){10}}
\linethickness{0.2mm}
\put(9,14){\oval(5,5)[l]}
\put(17,24){\line(-1,-1){8}}
\put(20,6){\oval(6,6)[b]}
\put(17,6){\line(-1,1){6}}
\put(23,6){\line(1,2){2}}
\put(20,24){\oval(6,6)[t]}
\put(17,24){\circle{2}}
\put(23,24){\line(0,-1){4}}
\put(28,10){\line(-1,2){5}}
\put(28,0){\line(0,1){10}}
\put(26,17){\line(0,1){13}}

\put(30,15){\mbox{$=$}}
\linethickness{0.4mm}
\put(43,0){\line(0,1){9}}
\put(43,14){\line(0,1){16}}
\linethickness{0.2mm}
\put(36,16){\oval(4,4)[t]}
\put(36,13){\oval(4,4)[b]}
\put(34,13){\line(0,1){4}}
\put(38,16){\line(1,-1){10}}
\put(48,6){\line(0,-1){6}}
\put(42,25){\line(-1,-3){3}}
\put(41,21){\circle{2}}
\put(44,27){\line(1,1){4}}

\put(50,15){\mbox{$=$}}
\linethickness{0.4mm}
\put(60,0){\line(0,1){16}}
\put(60,21){\line(0,1){9}}
\linethickness{0.2mm}
\put(61,30){\oval(8,8)[br]}
\put(59,22){\oval(5,8)[l]}
\put(59,18){\line(1,0){5}}
\put(64,14){\oval(8,8,)[tr]}
\put(68,14){\line(0,-1){14}}

\put(80,25){\mbox{\small(\ref{satz511}):}}
\linethickness{0.4mm}
\put(90,25){\line(0,1){5}}
\put(90,21){\line(0,-1){8}}
\put(90,0){\line(0,1){10}}
\linethickness{0.2mm}
\put(94,26){\line(0,1){4}}
\put(91,26){\oval(6,6)[br]}
\put(88,21){\oval(4,4)[l]}
\put(88,23){\line(1,0){4}}
\put(91,16){\oval(6,6)[tr]}
\put(99,16){\line(0,1){14}}
\put(99,16){\line(-2,-1){13}}
\put(97,4){\oval(6,6)[b]}
\put(100,4){\line(-1,2){4}}
\put(91,4){\oval(6,6)[tr]}
\put(89,10){\oval(6,6)[bl]}

\put(104,15){\mbox{$=$}}
\linethickness{0.4mm}
\put(112,0){\line(0,1){30}}
\linethickness{0.2mm}
\put(118,4){\line(1,2){3}}
\put(126,20){\line(-1,-2){3}}
\put(126,4){\line(-1,2){8}}

\put(118,20){\line(0,1){10}}
\put(126,20){\line(0,1){10}}

\put(122,4){\oval(8,8)[b]}

\end{picture}
\caption{\label{bildsatz5} 
Two of the relations of $\BRib$. The small circle denotes a ribbon twist
$\twist'$.}

\end{figure}

The proof is a simple verification. Some of the pictorial calculations are given
in figure \ref{bildsatz5}. Because of (\ref{satz53}) and
(\ref{satz54}) only $\tau^{\downarrow+},\tau^{\uparrow+}$
are needed as generators. This reduces the numbers of relations because
(\ref{satz53}) and (\ref{satz54}) turn (\ref{satz59}) and (\ref{satz510})
into identities that involve only $\ARib$ operations. 

\begin{figure}[ht]
\begin{picture}(150,40)

\put(20,20){\circle*{2}}
\linethickness{0.1mm}
\put(20,20){\line(-1,0){20}}
\linethickness{0.2mm}
\put(5,17){\line(0,1){6}}
\put(10,17){\line(0,1){6}}
\put(15,17){\line(0,1){6}}
\linethickness{0.4mm}
\put(3,23){\line(0,1){10}}
\put(3,23){\line(1,0){25}}
\put(3,17){\line(0,-1){10}}
\put(3,17){\line(1,0){25}}
\put(3,7){\line(1,0){40}}
\put(3,33){\line(1,0){40}}
\put(43,7){\line(0,1){26}}
\put(28,17){\line(0,1){6}}
\put(15,25){\mbox{{\small Schlingel}}}
\put(0,10){\mbox{$\downarrow$}}
\put(0,30){\mbox{$\uparrow$}}

\put(70,20){\mbox{$\mapsto$}}

\put(100,20){\circle*{2}}
\linethickness{0.2mm}
\put(100,20){\oval(20,20)[l]}
\put(100,20){\oval(10,10)[l]}
\put(100,20){\oval(30,30)[l]}
\put(100,35){\line(1,0){5}}
\put(100,30){\line(1,0){5}}
\put(100,25){\line(1,0){5}}
\put(100,15){\line(1,0){5}}
\put(100,10){\line(1,0){5}}
\put(100,5){\line(1,0){5}}

\linethickness{0.4mm}
\put(105,3){\line(0,1){14}}
\put(105,23){\line(0,1){14}}
\put(105,3){\line(1,0){20}}
\put(105,37){\line(1,0){20}}
\put(105,17){\line(1,0){5}}
\put(105,23){\line(1,0){5}}
\put(110,17){\line(0,1){6}}
\put(125,3){\line(0,1){34}}

\end{picture}
\caption{\label{rpbild} Deforming a horizontal diagram:
On the left the original diagram with a radar beam and arrows indicating the direction of
deformation. The result is shown on  the right.}
\end{figure}
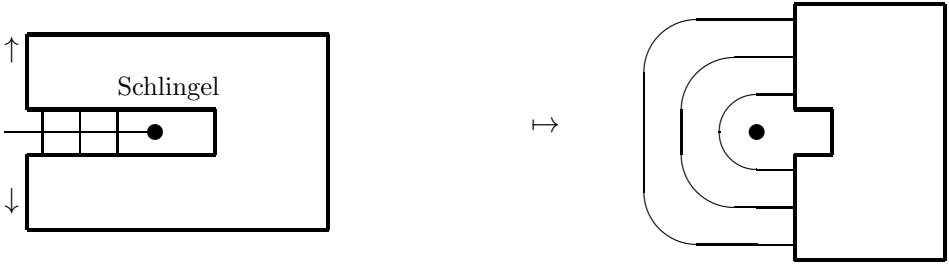

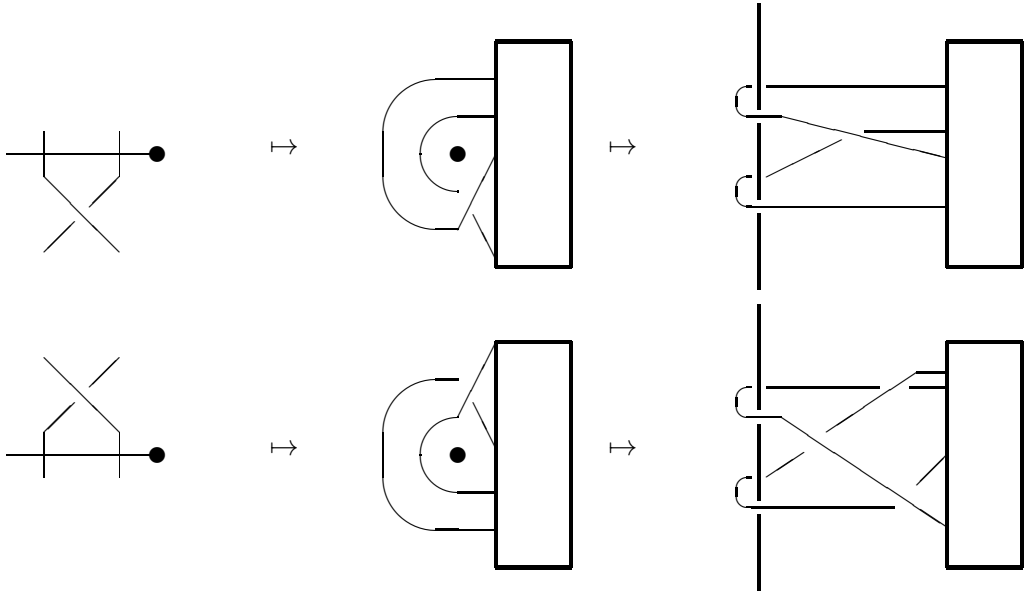
\begin{figure}[ht]
\begin{picture}(150,80)

\put(20,60){\circle*{2}}
\put(20,20){\circle*{2}}
\linethickness{0.1mm}
\put(0,20){\line(1,0){20}}
\put(0,60){\line(1,0){20}}
\linethickness{0.1mm}
\put(5,63){\line(0,-1){6}}
\put(15,63){\line(0,-1){6}}
\put(5,57){\line(1,-1){10}}
\put(5,47){\line(1,1){4}}
\put(15,57){\line(-1,-1){4}}

\put(5,23){\line(0,-1){6}}
\put(15,23){\line(0,-1){6}}
\put(5,33){\line(1,-1){10}}
\put(5,23){\line(1,1){4}}
\put(15,33){\line(-1,-1){4}}

\put(35,60){\mbox{$\mapsto$}}
\put(35,20){\mbox{$\mapsto$}}

\put(60,60){\circle*{2}}
\put(60,20){\circle*{2}}
\linethickness{0.2mm}
\put(60,60){\oval(10,10)[l]}
\put(60,60){\oval(20,20)[l]}

\put(60,65){\line(1,0){5}}
\put(60,70){\line(1,0){5}}

\put(60,15){\line(1,0){5}}
\put(60,10){\line(1,0){5}}

\put(60,50){\line(1,2){5}}
\put(65,46){\line(-1,2){3}}

\put(60,25){\line(1,2){5}}
\put(65,21){\line(-1,2){3}}

\put(60,20){\oval(10,10)[l]}
\put(60,20){\oval(20,20)[l]}

\linethickness{0.4mm}
\put(65,45){\line(0,1){30}}
\put(75,45){\line(0,1){30}}
\put(65,45){\line(1,0){10}}
\put(65,75){\line(1,0){10}}

\put(65,5){\line(0,1){30}}
\put(75,5){\line(0,1){30}}
\put(65,5){\line(1,0){10}}
\put(65,35){\line(1,0){10}}

\put(80,60){\mbox{$\mapsto$}}
\put(80,20){\mbox{$\mapsto$}}

\linethickness{0.4mm} 

\put(100,80){\line(0,-1){14}}
\put(100,64){\line(0,-1){10}}
\put(100,52){\line(0,-1){10}}

\put(100,40){\line(0,-1){14}}
\put(100,24){\line(0,-1){10}}
\put(100,12){\line(0,-1){10}}

\put(125,45){\line(0,1){30}}
\put(135,45){\line(0,1){30}}
\put(125,45){\line(1,0){10}}
\put(125,75){\line(1,0){10}}

\put(125,5){\line(0,1){30}}
\put(135,5){\line(0,1){30}}
\put(125,5){\line(1,0){10}}
\put(125,35){\line(1,0){10}}

\linethickness{0.2mm} 

\put(99,67){\oval(4,4)[l]}
\put(99,55){\oval(4,4)[l]}
\put(99,27){\oval(4,4)[l]}
\put(99,15){\oval(4,4)[l]}
\put(99,65){\line(1,0){4}}
\put(99,53){\line(1,0){4}}
\put(99,25){\line(1,0){4}}
\put(99,13){\line(1,0){4}}

\put(101,69){\line(1,0){24}}
\put(101,53){\line(1,0){24}}
\put(103,65){\line(4,-1){22}}
\put(101,57){\line(2,1){10}}
\put(114,63){\line(1,0){11}}

\put(103,25){\line(3,-2){22}}
\put(101,29){\line(1,0){15}}
\put(120,29){\line(1,0){5}}
\put(121,31){\line(-3,-2){12}}
\put(121,31){\line(1,0){4}}

\put(103,13){\line(1,0){15}}
\put(125,20){\line(-1,-1){4}}

\put(101,17){\line(3,2){5}}

\end{picture}
\caption{\label{vzr} The pictures in the upper and lower row 
differ only by the position of a single crossing relative to the chosen radar beam.
Irrelevant parts of the diagram are omitted.  
The second mapping associates a standard diagram to the horizontal diagram.
 }
\end{figure}

\begin{satz}
The set  ${\cal F}:=\{\tau^{\downarrow+},\tau^{\uparrow+}\}$ generates
the action pair $(\BRib,\ARib)$ with relations
(\ref{satz55})-(\ref{satz58}), (\ref{satz511})-(\ref{satz514}).                                                                                                                                                                                                                                                                                                                                                                                                                                                                                                                                                                                                                                                                                                                        
\end{satz}
\begin{bew}
Tangles in the cylinder may be interpreted as ordinary tangles with a fixed additional
strand. The question of equivalence of diagrams can thus be reduced to the
situation in $\RR^3$ \cite{tu}. However, ordinary Markov moves may easily produce
diagrams that are no longer products of our generators.
We therfore need a method to produce a standard form (a product of generators)
from an arbitrary diagram.
There are several such methods. We use the R-process. It is based on 
horizontal diagrams, i.e. regular projections of cylinder links on a horizontal 
plane. In contrast we call the diagrams used sofar standard diagrams.
In  a horizontal diagram the cylinder axis is just a point. To avoid upper and lower
end points from being projected on the same point we shift them in 
opposite directions parallel to the second coordinate axis.
Upon multiplication we have to join such endings with horizontal ribbons. Figure
\ref{neupro} displays an example. 

\unitlength1.7mm
\begin{figure}[ht]
\begin{picture}(70,25)
\linethickness{0.4mm}
\put(5,5){\line(0,1){8}}
\put(5,15){\line(0,1){10}}
\linethickness{0.2mm}
\put(4,16){\oval(4,4)[l]}
\put(4,14){\line(1,0){3}}
\put(7,12){\oval(4,4)[tr]}
\put(9,12){\line(0,-1){7}}
\put(8,20){\line(0,1){5}}
\put(6,20){\oval(4,4)[br]}
\put(11,25){\line(1,-4){5}}
\put(11,5){\line(1,4){2}}
\put(16,25){\line(-1,-4){2}}

\put(20,15){\mbox{$\leftrightarrow$}}

\put(30,15){\circle*{2}}
\put(30,15){\oval(6,6)[l]}
\put(30,21){\oval(6,6)[br]}
\put(30,9){\oval(6,6)[tr]}
\put(35,21){\line(1,-3){4}}
\put(35,9){\line(1,3){1.5}}
\put(39,21){\line(-1,-3){1.5}}

\put(45,20){\mbox{\small up}}
\put(45,10){\mbox{\small down}}

\end{picture}
\caption{\label{neupro} A simple  example of a horizontal diagram }
\end{figure}
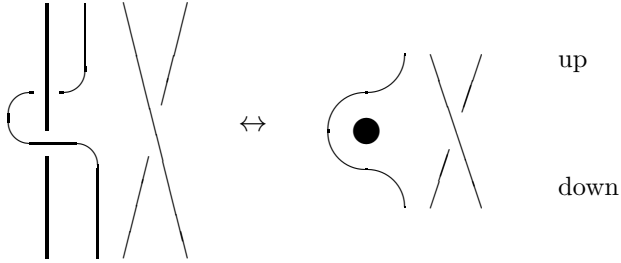

Let such a horizontal diagram be given and 
choose a line (the radar beam) from the point of the axis and 
extending into the left half plane such that it hits the ribbons transversal
and avoids  crossings. The tangle is then deformed away from the beam 
until all of its nontrivial part is located in the right half plane
as indicated in figure \ref{rpbild}.
We may assign a standard diagram to the result by drawing a sequence of $\tau$ morphisms
for every circle surrounding the axis such that the innermost circle corresponds to the
lowest $\tau$. Diagrams that differ by some kind of Reidemeister move that
takes place either above or below the radar beam are transformed to 
standard diagrams that are related by precisely the same move at a different position.
It remains to discuss how the result depends on the choice of the
radar beam. Essentially, there are only two relevant possibilities that
correspond to Reidemeister moves of types II and III. We concentrate on the
type III move.  Consider two situations differing only by the position of a single 
crossing with respect to the radar beam. 
The beam may either be above or below the crossing.
Figure  \ref{vzr} shows diagrams of these situations.
We demand that the tangle in the right half space is concentrated in a diagonal box
so that the connection points  with the axis surrounding circles are
projected both horizontally and vertically to the same order.
Then it is easy to determine how the parts fit together.
Comparing diagrams on the right of figure  \ref{vzr}
yields the four braid relation. A similar argument using a minimum (maximum)
lying above or below the radar beam yields the extremum twist relation.
Putting in orientations these relations are precisely
(\ref{satz55})-(\ref{satz58}), (\ref{satz511})-(\ref{satz514}).                                                                                                                                                                                                                                                                                                                                                                                                                                                                                                                                                                                                                                                                                                                        
\end{bew}

\unitlength1mm

\unitlength1mm

\begin{satz}
There is a unique tensor functor between strict action pairs
$F:(\BRib_\akat,\ARib_\akat)\rightarrow(\bkat,\akat)$ such that
$F_\akat$ is Turaev's functor, the  functorial isomorphism $\omega$ 
is trivial	and one has
\begin{eqnarray}
F_\bkat(\tau_X^{\downarrow\pm})&=&t_X^{\pm1}\\
F_\bkat(\tau_X^{\uparrow\pm})&=&t_{X^\ast}^{\pm1}
\end{eqnarray}
\end{satz}
\begin{bew}
Uniqueness is clear because $F_\bkat$ if fixed on generators.
To prove existence one has to check compatibility with the relations
given above. This is done by straightforward graphical computations
which are however too long to be displayed here.
\end{bew}

\section{Cylinder braided action pairs with Points}

Until now we have no possibility to represent the diagram $e_0$
of figure \ref{agenerat} which plays a crucial role in the study of some
B-type knot algebras. 
The point structure	discussed in this section fills the gap.

\begin{de}
A point structure on a CBTC $(\bkat,\akat)$ (where $\akat$ is rigid)
consists of a point morphism
$\bnull_X\in\Mor_\bkat(1,X)$ and copoint morphisms 
$\dnull_X\in\Mor_\bkat(X,1)$ such that the following axioms are fulfilled. 
\begin{eqnarray}
d^0_Yf&=&d^0_X\quad fb^0_X=b^0_Y\qquad\forall f\in{\rm Mor}_\akat(X,Y)\\
d_X(\bnull_{X^\ast}\otimes\id_X)&=&\dnull_X\\
(\dnull_X\otimes\id_{X^\ast})b_X&=&\bnull_{X^\ast}\\
\bnull_{X\otimes Y}&=&(\bnull_X\otimes\id_Y)\bnull_Y\\
\dnull_{X\otimes Y}&=&\dnull_Y(\dnull_X\otimes\id_Y)\\
\bnull_X&=&t_X\bnull_X\label{aufdreh}\\
\dnull_X&=&\dnull_Xt_X\\
(t_Y\otimes\id_X)c_{X,Y}(\bnull_X\otimes\id_Y)&=&
c^{-1}_{Y,X}(\bnull_X\otimes\id_Y)t_Y\\
t_Y(\dnull_X\otimes\id_Y)&=&\dnull_X c_{Y,X}(t_Y\otimes\id_X)c_{X,Y}
\end{eqnarray}
\end{de}
Some simple consequences are:
\begin{eqnarray}
\dnull_{X^\ast}&=&d_Xc_{X,X^\ast}(\theta_X\bnull_X\otimes\id_{X^\ast}) \\
\bnull_{X^\ast}&=&(\dnull_{X^\ast}\otimes\id_X)
(\id_{X^\ast}\otimes\theta^{-1}_X)
c^{-1}_{X\ast,X}b_X
\end{eqnarray}

A point structure is the B-type analog of duality (rigidity).

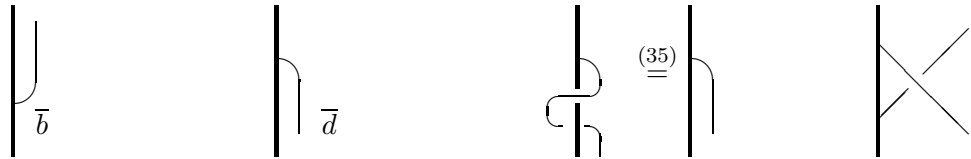
\begin{figure}[ht]
\begin{picture}(150,20)

\linethickness{0.4mm}
\put(5,0){\line(0,1){20}}
\linethickness{0.2mm}
\put(5,10){\oval(6,6)[br]}
\put(8,10){\line(0,1){8}}
\put(8,3){\mbox{$\bnull$}}

\linethickness{0.4mm}
\put(40,0){\line(0,1){20}}
\linethickness{0.2mm}
\put(40,10){\oval(6,6)[tr]}
\put(43,10){\line(0,-1){7}}
\put(46,3){\mbox{$\dnull$}}

\linethickness{0.4mm}
\put(80,0){\line(0,1){7}}
\put(80,20){\line(0,-1){11}}
\linethickness{0.2mm}
\put(80,10){\oval(6,6)[tr]}
\put(81,10){\oval(4,4)[br]}
\put(78,6){\oval(4,4)[l]}
\put(78,8){\line(1,0){3}}
\put(83,3){\line(0,-1){3}}
\put(81,2){\oval(4,4)[tr]}
\put(88,10){\mbox{$\stackrel{(\ref{aufdreh})}{=}$}}

\linethickness{0.4mm}
\put(95,0){\line(0,1){20}}
\linethickness{0.2mm}
\put(95,10){\oval(6,6)[tr]}
\put(98,10){\line(0,-1){7}}

\linethickness{0.4mm}
\put(120,0){\line(0,1){20}}
\linethickness{0.2mm}
\put(120,15){\line(1,-1){12}}
\put(120,5){\line(1,1){4}}
\put(126,11){\line(1,1){6}}

\end{picture}
\caption{\label{pbild} 
Point and copoint of  $P\BRib$}

\end{figure}

$(\BRib,\ARib)$ has no point structure.
We define $(P\BRib,\ARib)$ as an extension where
ribbons are allowed to end at the cylinder axis.
Figure \ref{pbild} display the point and copoint morphisms.
Note that points (i.e endings of ribbons on the axis) do not commute, i.e there
is no way to simplify the picture on the right of figure \ref{pbild}.

\unitlength1mm

\section{Skein relations}

In  $P\BRib$ one can impose skein relations that generalise those of
the Kauffman polynomial:

\begin{eqnarray}
c-c^{-1}&=&\delta(1-bd)\label{SK1}\\
cb&=&\lambda b\qquad dc=\lambda d \label{SK2}\\
db&=& A_0\label{SK3}\\
t^{-1}&=&\alpha t+\beta+\gamma b^0d^0 \label{SK4}\\
d^0b^0&=&x_0 \label{SK5}\\
db^0d^0b&=&x'_0 \label{SK6}\\
d(t\otimes\id)b&=&A_1 \label{SK7}\\
d(t^{-1}\otimes\id)b&=&A_{-1} \label{SK8}\\
(d^0\otimes\id)c(b^0\otimes\id)&=& \epsilon+\mu t+\nu b^0d^0\label{SK9}
\end{eqnarray}

The parameters are 
$\delta,A_0,\lambda,\alpha,\beta,\gamma,A_1,A_{-1},x_0,x'_0,\epsilon,\mu,\nu$.

Assuming that the annihilator ideals of the generators vanish we can 
derive a set of relations between these parameters.
As in the case of the A-type category of the usual Kauffman polynomial one has
\[A_0\delta-\delta=\lambda-\lambda^{-1}\]
We have 
$d^0=d^0t^{-1}=\alpha d^0+\beta d^0+\gamma d^0b^0d^0=(\alpha+\beta+\gamma x_0)d^0$.
and hence: \[1=\alpha+\beta+\gamma x_0\]
Similarly:
\[A_{-1}=\alpha A_1+\beta A_0+\gamma x_0'\]
Multiplying  $\lambda^{-1}(t^{-1}\otimes\id)b=c(t\otimes\id)b$
with $d$ we obtain \[A_1\lambda^2=A_{-1}\]
Next, we calculate  $\gamma x_0b^0d^0=\gamma b^0d^0b^0d^0=b^0d^0(t^{-1}-\alpha t-\beta)=
b^0d^0(1-\alpha-\beta)$ and obtain
\[\gamma x_0=1-\alpha-\beta\]
Similarly $x'_0=db^0d^0b=\gamma^{-1}(d(t^{-1}\otimes1)b-
\alpha d(t\otimes 1)b-\beta db)$:
\[\gamma x'_0=A_{-1}-\alpha A_1-\beta A_0\]
Finally, tensor  (\ref{SK4})  with $c$ and multiply with  $d\otimes\id$
from the left and with  $b\otimes\id$ from the right. The result may be brought to 
a form which resembles  (\ref{SK9}). Comparing coefficients one obtains:
\begin{eqnarray*}
\nu&=&-\alpha\lambda\\
\mu&=&\gamma^{-1}(\alpha\delta-\alpha^2\lambda+\lambda^{-1})\\
\epsilon&=&-\gamma^{-1}(\alpha\beta\lambda+\alpha\delta A_1+\beta\lambda^{-1})
\end{eqnarray*}

Only 4 of  13 parameters survive. We may reduce this number even more, if we demand
that $x_0=x'_0$.

Link in the solid torus are endomorphisms of the 
 $0$-object.  Kauffman's Theory \cite{kaufpoly} can be used to eliminate ordinary
 braidings $c$. The remaining tangles can be simplified using 
  (\ref{SK9}). Therefore the skein relations suffice to calculate a 
  cylinder generalisation of Kauffman´s polynomial.

Just as with Kauffman's original polynomial this link invariant may also be obtained
as a writhe normalisation of a Markov trace on a B-type generalisation of
Birman-Wenzl algebra \cite{rhofull}.

For special values of the parameters it may also be derived from tensor
representations of $P\BRib$. Upto now only two nontrivial tensor representations 
have been  found. They are based on  tensor representations of the 
B-type braid group that use R-matrices of the orthogonal quantum groups 
(found by tom Dieck \cite{tD2} for ${\cal U}_q{\rm so}_3$
and by myself for ${\cal U}_q{\rm so}_5$ (unpublished)).

\end{document}